\documentclass[leqno,draft]{article}

\newtheorem{theorem}{Theorem}
\newtheorem{lemma}[theorem]{Lemma}
\newtheorem{proposition}[theorem]{Proposition}
\newtheorem{definition}[theorem]{Definition}
\newtheorem{corollary}[theorem]{Corollary}

\newcommand{\begintheorem}{\addtocounter{equation}{1}\begin{theorem}}
\newcommand{\beginlemma}{\addtocounter{equation}{1}\begin{lemma}}
\newcommand{\beginproposition}{\addtocounter{equation}{1}\begin{proposition}}
\newcommand{\begindefinition}{\addtocounter{equation}{1}\begin{definition}}
\newcommand{\begincorollary}{\addtocounter{equation}{1}\begin{corollary}}

\begin{document}

\title{Elements of harmonic analysis, 3}

\author{Stephen William Semmes	\\
	Rice University		\\
	Houston, Texas}

\date{}

\maketitle

        These informal notes are based on a course given at Rice
University in the spring semester of 2004, and much more information
can be found in the references.

\section*{Locally compact abelian topological groups}

	Let $A$ be an abelian group.  Thus $A$ is a set equipped with
a binary operation $+$ which is commutative and associative, there
is an identity element $0 \in A$ such that $0 + a = a$ for all $a \in A$,
and each $a \in A$ has an inverse $-a$ characterized by $a + (-a) = 0$.
As basic examples, the integers, real numbers, and complex numbers
are abelian groups under addition, and for each positive integer $n$
we have the integers modulo $n$, a cyclic group with $n$ elements.

	Let us also assume that $A$ is a topological space, which is
to say that certain subsets of $A$ are designated as open subsets.  As
usual one requires that the empty set and $A$ itself are open subsets
of $A$, that the intersection of any finite collection of open subsets
is an open subset, and that the union of any collection of open
subsets is an open subset.  Once the open subsets are selected, the
closed subsets are defined to be the complements in $A$ of the open
subsets.  Various standard notions, such as continuity at a point
of a mapping between two topological spaces, can be defined in terms
of the open subsets through standard methods.

	To say that $A$ is a topological group means that the group
structure and topology are compatible in a natural way.  Specifically
the group operation $+$ should be continuous as a mapping from $A
\times A$ into $A$, and $a \mapsto -a$ should be continuous as a
mapping from $A$ to itself.  This implicitly uses the product topology
on $A \times A$ induced by the given topology on $A$, in which a
subset of $A \times A$ is open if it is the union of products of open
subsets of $A$.  It is customary to require that $A$ be a Hausdorff
topological space, which is equivalent in this setting to the
requirement that $\{0\}$ be a closed subset of $A$.

	In any topological space a subset $K$ is said to be compact if
every open covering of $K$ in the space admits a finite subcovering,
i.e., if for every family $\{U_\iota\}_{\iota \in I}$ of open subsets
of the topological space such that
\begin{equation}
	K \subseteq \bigcup_{\iota \in I} U_\iota
\end{equation}
there is a finite collection $\iota_1, \ldots, \iota_l$ of indices in
$I$ such that
\begin{equation}
	K \subseteq U_{\iota_1} \cup \cdots \cup U_{\iota_l}.
\end{equation}
A topological space is said to be locally compact if for each
point $x$ in the space there is an open subset $W$ and a compact
subset $K$ of the space such that
\begin{equation}
	x \in W \subseteq K.
\end{equation}
A locally compact abelian topological group is an abelian topological
group which is locally compact as a topological space.  Of course
local compactness at the identity element $0$ implies local
compactness at every point because group translations define
homeomorphisms.

	As in \cite{3}, it is simpler to say ``LCA group'' in place of
locally compact abelian topological group.  The integers and the
integers modulo $n$ are natural examples of LCA groups equipped with
their discrete topologies, in which every subset is considered to be
open.  For the real and complex numbers one can use their standard
topologies, induced by the usual Euclidean metrics, to get LCA groups.
One can also consider the nonzero complex numbers using multiplication
as the group operation and the usual topology.  If one takes the
complex numbers with modulus $1$ using multiplication as the group
operation and the usual topology, one gets a compact LCA group.

	Fix an integer $n \ge 2$, and consider the group consisting of
sequences $x = \{x_j\}_{j=1}^\infty$ such that each $x_j$ is an
integer modulo $n$.  One might as well say that each $x_j$ is an
integer such that $0 \le x_j \le n-1$, and the sum of two elements
$x$, $y$ in the group is defined by adding each term modulo $n$.  If
$x$ is an element of the group and $l$ is a positive integer, then the
$l$th standard neighborhood around $x$ is defined to be the set of $y$
in the group such that $x_j = y_j$ when $1 \le j \le l$.  This leads
to a topology on the space in which a subset of the group is open if
for each point $x$ in the subset there is a standard neighborhood
around $x$ which is contained in the subset.  Well known results
in topology imply that this space is compact with respect to this
topology, and in fact it is homeomorphic to the Cantor set.

	It is easy to see that this example defines a compact LCA
group.  Namely, the group operations are continuous with respect to
the topology just defined.  This is a nice example where the
topological dimension is equal to $0$, which is to say that the space
is totally disconnected, with no connected subsets with at least two
elements.  At the same time the topology is not the discrete topology.

	Let $A$ be a LCA group.  A basic object of interest associated
to $A$ is a translation-invariant integral, which is a linear mapping
from the vector space of complex-valued continuous functions $f(x)$ on
$A$ with compact support into the complex numbers such that the
integral of $f(x + a)$ is equal to the integral of $f(x)$ for all $a
\in A$, the integral of a real-valued function is a real number, the
integral of a nonnegative real-valued function is a nonnegative real
number, and the integral of $f$ is positive if $f$ is a nonnegative
real-valued continuous function on $A$ with compact support such that
$f(x) > 0$ for some $x \in A$.  In the examples described earlier such
an integral can be defined explicitly, in terms of sums, classical
Riemann integrals, or simple generalizations of Riemann integrals for
the spaces of sequences modulo $n$.  A general theorem states that any
LCA group $A$ has such an invariant integral, and that this integral
is unique except for multiplying it by a positive real number.

	Let $A$ be an LCA group.  By a character on $A$ we mean a
continuous group homomorphism from $A$ into the group of complex
numbers with modulus $1$ with respect to multiplication.  Sometimes
one may wish to consider unbounded characters more generally, which
are continuous homomorphisms from $A$ into the nonzero complex numbers
with respect to multiplication.  Note that any bounded subgroup of the
nonzero complex numbers with respect to multiplication is contained in
the complex numbers with modulus $1$, as one can easily verify.
Thus a bounded continuous homomorphism from $A$ into the group of
nonzero complex numbers is a character, and in particular every
continuous homomorphism from $A$ into the nonzero complex numbers
is a character when $A$ is compact.

	If $f(x)$ is a complex-valued continuous function on $A$ with
compact support, or an integrable function more generally, one can
define its Fourier transform $\widehat{f}(\phi)$ by saying that if
$\phi$ is a character on $A$, then $\widehat{f}(\phi)$ is the integral
of $f$ times the complex conjugate of $\phi$, using a fixed invariant
integral on $A$ as discussed previously.  If $A$ is not compact, then
one can extend this to a Fourier--Laplace transform by allowing
unbounded characters, at least when $f$ has compact support or
sufficient integrability properties.  For bounded characters one has
the usual inequality which states that $|\widehat{f}(\phi)|$ is less than
or equal to the integral of $|f|$.

	Many classical aspects of Fourier analysis work in this
setting.  A basic point is that the Fourier transform diagonalizes
translation operators, which means that if $a \in A$ and $f(x)$
is a continuous function on the group with compact support, or an
integrable function on the group, then the Fourier transform of
$f(x - a)$ at the character $\phi$ is equal to $\overline{\phi(a)}$
times the Fourier transform of $f(x)$ at $\phi$.  One can also define
convolutions in the usual way, using the invariant integral on $A$,
and the Fourier transform of a convolution is equal to the
product of the corresponding Fourier transforms.


\begin{thebibliography}{4}



\bibitem {1} W.~Hurewicz and H.~Wallman, {\it Dimension Theory},
Princeton Mathematical Series {\bf 4}, Princeton University Press,
1941.

\bibitem {2} E.~Hewitt and K.~Ross, {\it Abstract Harmonic Analysis I},
second edition, Springer-Verlag, 1979.

\bibitem {3} W.~Rudin, {\it Fourier Analysis on Groups}, Wiley, 1962.

\bibitem {4} E.~Stein and G.~Weiss, {\it Introduction to Fourier
Analysis on Euclidean Spaces}, Princeton Mathematical Series {\bf 32},
Princeton University Press, 1971.






\end{thebibliography}
\end{document}